\documentclass[10pt]{article}
\usepackage{amssymb}
\usepackage{amsmath,amssymb}\usepackage{cite}\usepackage{mathrsfs}\usepackage[amsmath,thmmarks]{ntheorem}
\usepackage{listings}
\usepackage[titletoc]{appendix}\allowdisplaybreaks

\makeatletter
\renewcommand\theequation{\thesection.\@arabic\c@equation}
\makeatother

\newtheorem{thm}{ Theorem}[section]%
\newtheorem{lem}[thm]{ Lemma}%
\newtheorem{cor}[thm]{ Corollary}%
\newtheorem{Que}[thm]{Question}%
\input cyracc.def

\def\f{\noindent}

\def\demo{\f{\bf Proof}\hskip10pt}

\def\qed{\hfill $\Box$}

\begin{document}
\title{\textbf{A Note On The Cross-Sperner Families}}
\author{Junyao Pan\,
 \\\\
School of Sciences, Wuxi University, Wu'xi, Jiangsu,\\ 214105
 People's Republic of China \\}
\date {} \maketitle

\baselineskip=16pt

\vskip0.5cm

{\bf Abstract:} Let $(\mathcal{F},\mathcal{G})$ be a pair of families of $[n]$, where $[n]=\{1,2,...,n\}$. If $A\not\subset B$ and $B\not\subset A$ hold for all $A\in\mathcal{F}$ and $B\in\mathcal{G}$, then $(\mathcal{F},\mathcal{G})$ is called a Cross-Sperner pair. P. Frankl and Jian Wang introduced the extremal problem that $m(n)={\rm{max}}\{|\mathcal{I}(\mathcal{F},\mathcal{G})|:\mathcal{F},\mathcal{G}\subset2^{[n]}~{\rm{are~cross}}$-${\rm{sperner}}\}$, where $\mathcal{I}(\mathcal{F},\mathcal{G})=\{A\cap B:A\in\mathcal{F},B\in\mathcal{G}\}$. In this note, we prove that $m(n)=2^n-2^{\lfloor\frac{n}{2}\rfloor}-2^{\lceil\frac{n}{2}\rceil}+1$ for all $n>1$. This solves an open problem proposed by P. Frankl and Jian Wang.



\section {Introduction}

Let $[n]=\{1,2,...,n\}$ denote the standard $n$-element set and $2^{[n]}$ stand for the power set consisting of the $2^n$ subsets of $[n]$. A subset of $2^{[n]}$ is called a \emph{family}. The extremal problems of families of finite sets have always been an interesting research direction, the readers interested in this aspect can refer to \cite{FT}. A family $\mathcal{F}\in2^{[n]}$ is an antichain (also known as a Sperner family) if for all distinct $F, G \in \mathcal{F}$, neither $F\subseteq G$ nor $G \subseteq F$, i.e. $F$ and $G$ are incomparable. One of the principal results in extremal combinatorics is Sperner's theorem in \cite{S1}, which states that the largest size of an antichain in $2^{[n]}$ is $\dbinom{n}{\lfloor\frac{n}{2}\rfloor}$. It is natural to consider a generalisation of Sperner's theorem to multiple families of sets, such as a pair of families $(\mathcal{F}, \mathcal{G})$ is said to be Cross-Sperner if $A\not\subset B$ and $B\not\subset A$ hold for all $A\in\mathcal{F}$ and $B\in\mathcal{G}$. The research of such objects goes back to the 1970s when Seymour \cite{S} deduced from a result of Kleitman \cite{K1} that a Cross-Sperner pair $(\mathcal{F}, \mathcal{G})$ in $2^{[n]}$ satisfies $|\mathcal{F}|^{\frac{1}{2}}+|\mathcal{G}|^{\frac{1}{2}}\leq2^{\frac{1}{2}}$. In fact, there exist various extremal problems about Cross-Sperner families that have been studied, for example \cite{BK,BKMW,BMT,CM,E,GLPPS,K,K1,S,S1}.

Recently, P. Frankl and Jian Wang \cite{FW} by setting $\mathcal{I}(\mathcal{F},\mathcal{G})=\{A\cap B:A\in\mathcal{F},B\in\mathcal{G}\}$, they proposed the extremal problem that $m(n)={\rm{max}}\{|\mathcal{I}(\mathcal{F},\mathcal{G})|:\mathcal{F},\mathcal{G}\subset2^{[n]}~{\rm{are~Cross}}$-${\rm{Sperner}}\}$. Concurrently, they applied a skillful way to show that $m(n)=2^n-2^{d+1}+1$ for $n=2d$, for details see \cite[Theorem~5.5]{FW}. However, their approach is unable to determine $m(n)$ in case when $n$ is an odd number. Inspired by \cite[Example~5.4,Theorem~5.5]{FW}, they proposed the following problem.

\begin{Que}\label{pan1-1}\normalfont(\cite[Problem~5.6]{FW}~)
For $n=2d+1$, does $m(n)=2^n-2^{d+1}-2^d+1$ hold?
\end{Que}

See \cite[Example~5.4]{FW}, if $\mathcal{F}=\{A\cup Y:A\subsetneq X\}$ and $\mathcal{G}=\{B\cup X:B\subsetneq Y\}$ where $X\cup Y=[n]$ and $X\cap Y=\emptyset$, then we say that $(\mathcal{F},\mathcal{G})$ is a Cross-Sperner pair of \emph{type} $(X,Y)$. In this note, we prove the following theorem and thus we give an affirmative answer to Question\ \ref{pan1-1}.

\begin{thm}\label{pan1-6}\normalfont
Let $(\mathcal{F},\mathcal{G})$ be a Cross-Sperner pair of $[n]$ such that $|\mathcal{I}(\mathcal{F},\mathcal{G})|=m(n)$. Then $(\mathcal{F},\mathcal{G})$ is a Cross-Sperner pair of type $(X,Y)$ for some $X$ and $Y$ with $|X|=\lfloor\frac{n}{2}\rfloor$ and $|Y|=\lceil\frac{n}{2}\rceil$. In particular, $m(n)=2^n-2^{\lfloor\frac{n}{2}\rfloor}-2^{\lceil\frac{n}{2}\rceil}+1$ for all $n>1$.
\end{thm}

\section {Proof of Theorem 1.2}

Let $\dbinom{[n]}{k}$ denote the set of all $k$-subsets of $[n]$ for some $k$ with $1\leq k\leq n-1$. Next we start to prove Theorem\ \ref{pan1-6} from the following lemma.

\begin{lem}\label{pan2-1}\normalfont
Let $(\mathcal{F},\mathcal{G})$ be a Cross-Sperner pair of $[n]$ such that $|\mathcal{I}(\mathcal{F},\mathcal{G})|=m(n)$. Then, both $\mathcal{F}\cap\dbinom{[n]}{n-1}\neq\emptyset$ and $\mathcal{G}\cap\dbinom{[n]}{n-1}\neq\emptyset$ hold.
\end{lem}
\demo Proof by contradiction. Suppose that $\mathcal{F}\cap\dbinom{[n]}{n-1}=\emptyset$. Then $|A\cap B|\leq n-3$ for each $A\cap B\in\mathcal{I}(\mathcal{F},\mathcal{G})$ and thus $|\mathcal{I}(\mathcal{F},\mathcal{G})|\leq 2^{n-3}$. However, \cite[Example~5.4]{FW} indicate $|\mathcal{I}(\mathcal{F},\mathcal{G})|>2^{n-3}$, a contradiction. Therefore, $\mathcal{F}\cap\dbinom{[n]}{n-1}\neq\emptyset$. Likewise, $\mathcal{G}\cap\dbinom{[n]}{n-1}\neq\emptyset$.    \qed

\begin{cor}\label{pan2-2}\normalfont
Let $(\mathcal{F},\mathcal{G})$ be a Cross-Sperner pair of $[n]$ such that $|\mathcal{I}(\mathcal{F},\mathcal{G})|=m(n)$. Then, both $\bigcap\limits_{A\in\mathcal{F}}{A}\neq\emptyset$ and $\bigcap\limits_{B\in\mathcal{G}}{B}\neq\emptyset$ hold.
\end{cor}
\demo By Lemma\ \ref{pan2-1}, we may assume that there exist two $(n-1)$-subsets $A$ and $B$ of $[n]$ such that $A\in\mathcal{F}$ and $B\in\mathcal{G}$. According to the definition of Cross-Sperner, it follows that $i\in\bigcap\limits_{A\in\mathcal{F}}{A}$ and $j\in\bigcap\limits_{B\in\mathcal{G}}{B}$ where $\{i\}=[n]\setminus B$ and $\{j\}=[n]\setminus A$, as desired.    \qed

\begin{lem}\label{pan2-3}\normalfont
Let $(\mathcal{F},\mathcal{G})$ be a Cross-Sperner pair of $[n]$ such that $X=\bigcap\limits_{A\in\mathcal{F}}{A}$ and $Y=\bigcap\limits_{B\in\mathcal{G}}{B}$ and $X\cup Y=[n]$. Then $\mathcal{I}(\mathcal{F},\mathcal{G})\subseteq\mathcal{I}(\mathcal{F}',\mathcal{G}')$ and $|\mathcal{I}(\mathcal{F},\mathcal{G})|\leq |\mathcal{I}(\mathcal{F}',\mathcal{G}')|$ where $(\mathcal{F}',\mathcal{G}')$ is a Cross-Sperner pair of type $(X,Y)$.
\end{lem}
\demo Without loss of generality, we set $X=\bigcap\limits_{A'\in\mathcal{F}'}{A'}$ and $Y=\bigcap\limits_{B'\in\mathcal{G}'}{B'}$. By the definitions of Cross-Sperner and type $(X,Y)$, it follows that $\mathcal{F}\subseteq\mathcal{F}'$ and $\mathcal{G}\subseteq\mathcal{G}'$, as desired.    \qed

\begin{lem}\label{pan2-4}\normalfont
Let $(\mathcal{F},\mathcal{G})$ be a Cross-Sperner pair of $[n]$ such that $|\mathcal{I}(\mathcal{F},\mathcal{G})|=m(n)$. Suppose that $X=\bigcap\limits_{A\in\mathcal{F}}{A}$ and $Y=\bigcap\limits_{B\in\mathcal{G}}{B}$. Then, both $X\cap Y=\emptyset$ and $X\cup Y=[n]$ hold.
\end{lem}
\demo Assume that $m\in X\cap Y$. Set $\mathcal{F}'=\{A\setminus\{m\}:A\in\mathcal{F}\}$ and $\mathcal{G}'=\{B\setminus\{m\}:B\in\mathcal{G}\}$. Clearly, $|\mathcal{I}(\mathcal{F},\mathcal{G})|=|\mathcal{I}(\mathcal{F}',\mathcal{G}')|$. On the other hand, it is straightforward to see that $|\mathcal{I}(\mathcal{F}',\mathcal{G}')|\leq m(n-1)$. However, one easily checks that $m(n)>m(n-1)$, a contradiction. Therefore, $X\cap Y=\emptyset$ holds.

Pick a Cross-Sperner pair $(\mathcal{H},\mathcal{K})$ of $X\cup Y$ of type $(X,Y)$. Suppose that $X\cup Y\neq[n]$. Let $\Delta\subset[n]\setminus(X\cup Y)$. Define $\mathcal{H}_\Delta=\{A\cup\Delta:A\in\mathcal{H}\}$ and $\mathcal{K}_\Delta=\{B\cup\Delta:B\in\mathcal{K}\}$. It is obvious that $|\mathcal{I}(\mathcal{H}_\Delta,\mathcal{K}_\Delta)|=|\mathcal{I}(\mathcal{H},\mathcal{K})|$. Consider $\mathcal{F}$ and $\mathcal{G}$. Define $\mathcal{F}^\Delta=\{A\in\mathcal{F}:\Delta\subseteq A,A\subseteq(X\cup Y\cup\Delta)\}$ and $\mathcal{G}^\Delta=\{B\in\mathcal{G}:\Delta\subseteq B,B\subseteq(X\cup Y\cup\Delta)\}$.
It follows from Lemma\ \ref{pan2-3} that $$|\mathcal{I}(\mathcal{F}^\Delta,\mathcal{G}^\Delta)|\leq|\mathcal{I}(\mathcal{H}_\Delta,\mathcal{K}_\Delta)|=|\mathcal{I}(\mathcal{H},\mathcal{K})|.$$ On the other hand, for any two subsets $\Delta_1\subset[n]\setminus(X\cup Y)$ and $\Delta_2\subset[n]\setminus(X\cup Y)$, the Lemma\ \ref{pan2-3} implies that $\mathcal{I}(\mathcal{F}^{\Delta_1},\mathcal{G}^{\Delta_2})\subseteq\mathcal{I}(\mathcal{H}_{\Delta_1\cap\Delta_2},\mathcal{K}_{\Delta_1\cap\Delta_2})$. Hence, we deduce that $$|\mathcal{I}(\mathcal{F},\mathcal{G})|\leq\sum\limits_{\Delta\subset[n]\setminus(X\cup Y)}|\mathcal{I}(\mathcal{H}_\Delta,\mathcal{K}_\Delta)|= (2^{n-m}-1)|\mathcal{I}(\mathcal{H},\mathcal{K})|.$$ Since $|\mathcal{I}(\mathcal{H},\mathcal{K})|=2^m-2^{|X|}-2^{|Y|}+1$, it follows that $$|\mathcal{I}(\mathcal{F},\mathcal{G})|\leq (2^{n-m}-1)(2^m-2^{|X|}-2^{|Y|}+1)<2^n-2^{n-|X|}-2^{n-|Y|}+2^{n-m}.$$ Note that $2^n-2^{n-|X|}-2^{n-|Y|}+2^{n-m}=2^n-2^{n-|X|}-2^{|X|}+1-(2^{n-|Y|}-2^{|X|}-2^{n-m}+1)$. Since $n-|Y|>|X|$ and $n-|Y|>n-m$, we have $2^{n-|Y|}-2^{|X|}-2^{n-m}+1>0$. Therefore, $|\mathcal{I}(\mathcal{F},\mathcal{G})|<|\mathcal{I}(\mathcal{D},\mathcal{E})|$ where $\mathcal{D}$ and $\mathcal{E}$ are two Cross-Sperner families of $[n]$ of type $(D,E)$ with $|D|=|X|$ and $|E|=n-|X|$, a contradiction.   \qed

According to Lemma\ \ref{pan2-4} and the inequality $a+b\geq2\sqrt{ab}$ for $a\geq0$ and $b\geq0$, we obtain Theorem\ \ref{pan1-6} immediately.

\section{Acknowledgement}

We are very grateful to professor Stijn Cambie for correcting some mistakes about the literature related to corresponding problem in the initial version.

\end{document}